\documentclass[12pt,reqno]{article}
\usepackage{amsmath, amsthm, graphicx,amsfonts, amssymb,color}
\theoremstyle{definition}

\theoremstyle{remark}

\theoremstyle{example}

\numberwithin{equation}{section}
\usepackage{mathrsfs}
\newcommand{\scr}[1]{\mathscr #1}

\newcommand{\set}[1]{\left\{#1\right\}}
\newcommand{\R}{\mathbb R}
\newcommand{\eps}{\varepsilon}

\newcommand{\E}{\mathbb{E}}
\newcommand{\D}{\scr{D}}

\def\K{\scr K}

\renewcommand{\P}{\mathbb P}
\def\R{\mathbb R}
\def\Z{\mathbb Z}
\def\bg{\begin}
\def\be{\bg{equation}}
\def\de{\end{equation}}

\def\edar{\end{eqnarray}}

\def\lb{\label}
\def\ct{\cite}
\def\l{\left}\def\r{\right}
\def\fr{\frac}
\def\alp{\alpha}
\def\bt{\beta}

\def\dlt{\delta}

\def\eps{\epsilon}

\def\lmd{\lambda}
\def\Lmd{\Lambda}

\def\sgm{\sigma}

\def\fa{\forall}

\def\ex{\exists}

\def\ift{\infty}

\def\rar{\rightarrow}
\def\Rar{\Rightarrow}

\def\q{\quad}

\def\var{\text {\rm Var}}
\def\ess{{\rm ess}}

\def\lan{\langle}\def\ran{\rangle}

\def\[{\l[} \def\]{\r]}
\def\({\l(} \def\){\r)}
\def\|{\bigg|}

\def\hat{\widehat}

\newcommand{\rf}[1]{(\ref{#1})}

\def\var{{\rm{Var }}}

\def\P{\mathbb P}

 \def\beqlb{\begin{eqnarray}}\def\eeqlb{\end{eqnarray}}
 \def\beqnn{\begin{eqnarray*}}\def\eeqnn{\end{eqnarray*}}

\def\prodd{\displaystyle\prod}

\title{\bf  {Hitting Time Distributions for Denumerable Birth and Death Processes}
 \footnote{Research supported in part by Program for New Century Excellent Talents in University (NCET),
  973 Project(No 2006CB805901), NSFC(No 10721091)}}

\author{{\bf Yu Gong and  Yong-Hua Mao\footnote{Corresponding author: maoyh@bnu.edu.cn}}\\
{\small School of Mathematical Sciences, Beijing Normal University, }
\\
{\small Laboratory of Mathematics and Complex Systems, Ministry of Education}
\\
{\small Beijing 100875, People's Republic of China}
}

\date{}

\begin{document}

 \maketitle

\bg{abstract}  We proved the explicit formulas in Laplace transform of the hitting times for the birth and death processes on a denumerable
state space with $\ift$ the exit or entrance boundary. This extends the well known Keilson's theorem from finite state space to infinite state space.
 We also apply these formulas to the fastest strong stationary time for strongly ergodic
birth and death processes, and obtain the explicit convergence rate in separation.
 \end{abstract}



{\bf Running head}  Hitting time distribution for birth and death process

\section{Introduction}

In this paper, we will study the passage time between any two states of an irreducible birth and death process on the
nonnegative integers $\set{0,1,2,\cdots}$.
A well-known theorem states that the passage time from state $0$ to state $d(<\ift)$ is distributed as
a sum of $d$ independent exponential random variables with distinct rates. These rates are just the
non-zero eigenvalues of the associated generator for the process absorbed at state $d$.
This is a well-known theorem usually attributed to Keilson(\ct{ke}), and 
it may be traced back at least as far as Karlin and McGregor (\ct{km}).
See Diaconis and Miclo \ct{dm} for historical comments.

Very recently, Fill \ct{fill2}
gave a first stochastic proof for the result via the duality.
 An excellent application of this theorem is to the distribution of the fastest strong stationary time for an
 ergodic birth and death process on $\set{0,1,\cdots, d}$. 
And it is also the starting point of studying separation cut-off for birth and death processes in \ct{dsc}. 
By the similar method, Fill proved an analogue result for the
upward skip-free processes (\ct{fill3}).
Diaconis and Miclo \ct{dm} presented another probabilistic  proof for it, by using the ``differential operators'' for birth and death processes(\ct{fell}).

Consider a continuous-time birth and death process $(X_t)_{t\ge0}$ with generator $Q=(q_{ij})$ on $\Z_+$. The $(q_{ij})$ is as follows
\be\lb{eq0}
q_{ij}=\left\{
         \begin{array}{ll}
           b_i, & \hbox{for $j=i+1, i\ge0$;} \\
           a_i, & \hbox{for $j=i-1, i\ge1$;} \\
           -(a_i+b_i)& \hbox{for $j=i\ge1$;} \\
           -b_0& \hbox{for $j=i=0$;} \\
           0, & \hbox{for other $j\not=i$.}
         \end{array}
       \right.
\de
Here $a_i(i\ge1), b_i(i\ge0)$ be two sequences of positive numbers.

Let $T_{i,n}=\inf\set{t\ge0:X_t=n|X_0=i}$ be the hitting time of the state $n$ starting from the state $i$.
The well known theorem of Keilson(\ct{ke}) is the following.

\bg{thm}\lb{t01}
Let $\lmd_1^{(n)}<\cdots<\lmd_n^{(n)}$ be all (positive) $n$ eigenvalues of $-Q^{(n)}$, where
\be\lb{eq1}
Q^{(n)}=\left(
  \begin{array}{ccccccc}
    -b_0 & b_0 & 0 & 0 & \cdots & 0& 0 \\
    a_1 & -(a_1+b_1) & b_1 & 0 & \cdots & 0& 0 \\
    0 & a_2 & -(a_2+b_2) & b_2 & \cdots & 0& 0 \\
    \vdots & \ddots & \ddots & \ddots & \ddots & \ddots&  \vdots\\
    0 & 0 & 0& 0 & \cdots & a_{n-1} & -(a_{n-1}+{ b_{n-1}}) \\
  \end{array}
\right).
\de
Then $T_{0,n}$ is distributed as a sum of $n$ independent
exponential random variables with rate parameters $\set{\lmd_1^{(n)}, \cdots,\lmd_n^{(n)}}$. That is
\be\lb{e01}
\E e^{-sT_{0,n}}=\prod_{\nu=1}^n\fr{\lmd_\nu^{(n)}}{s+\lmd_\nu^{(n)}}, s\ge0.
\de
\end{thm}

 We will investigate the distribution of the hitting time $T_{i,n}$ for the birth and death process $X_t$ on the nonnegative integers. This includes four cases:
 \begin{description}
   \item[Case I:] $0\le i< n<\ift$;
  \item[Case II:] $0\le n<i\le N<\ift$, where $N$ is a reflecting state;
  \item[Case III:] $0\le i< n=\ift$;
  \item[Case IV:] $0\le n<i<\ift$.
 \end{description}

Cases I\&II are really easy consequence of Theorem \ref{t01} since the distributions are actually involved in finite states. This will be done in the
next section.
Indeed, for Case I, by using the property of the birth and death process and the strong Markov property, we can obtain the explicit
formula for any $0\le i< n<\ift$ from that of $T_{0,n}$ in Theorem \ref{t01}. See Corollary \ref{c01} below. For Case II, since $N$ is a reflecting
state, we can get the distribution of $T_{i,n}( 0\le n<i\le N)$ from Case I via turning left-side to right, that is, we can take the mapping on the state
space: $j\rar j': j'=N-j$. See Corollary \ref{c02} below in Section 2.

When we deal with the birth and death process on $\Z_+$, we will face the classification of the state $\ift$ at
infinity concerning uniqueness, due to Feller (\ct{fell}). See also \ct[Chapter 8]{and} for more details.
According to \ct{fell}, there are four types of the $\ift$ boundary:
regular, exit, entrance and natural boundaries.

 Define
\be\lb{mu}
\mu_0=1, \mu_i=\fr{b_0b_1\cdots b_{i-1}}{a_1a_2\cdots a_i}, i\ge 1,
\de
and $\mu=\sum_{i=1}^\ift\mu_i$.


Let us also define
\be\lb{rs}
R=\sum_{i=0}^\ift\fr{1}{\mu_ib_i}\sum_{j=0}^i\mu_j, S=\sum_{k=0}^\ift\fr{1}{\mu_kb_k}\sum_{i=k+1}^\ift\mu_i.
\de
The $\ift$ boundary is called exit if $R<\ift, S=\ift$; entrance if $R=\ift, S<\ift$.
 What we will do in this paper is to give distributions of $T_{i,\ift}$ for the birth and death process with the exit boundary and
 that of $T_{i,n} (i>n)$ for the entrance boundary.

 The another difficulty
 for infinite birth and death processes is obviously how about all the eigenvalues or the spectrum of the generator. By the spectral theory
 established in \ct{m04} and \ct{m06}, we can eventually overcome the difficulty. Briefly speaking, we give distributions of $T_{0,\ift}$(the life time) for
 the minimal birth and death process corresponding to $Q$ when $\ift$ is the exit boundary. We will use a procedure of approximation with $n\rar\ift$ to
 derive the distribution of $T_{0,\ift}$ from that of $T_{0,n}$ in Theorem \ref{t01}.
 To deal with the eigenvalues for birth and death processes in infinite state spaces, we should utilize the powerful theory of Dirichlet form.
Dirichlet form helps one obtain the variational formulas for eigenvalues, and more importantly provide the approximation procedure.
The similar situation appears when $\ift$ is the entrance boundary, from a view point of \ct{km} on the duality method. The duality method was used
successfully in \ct{cbk3} to study the estimation of the principal eigenvalue for birth and death processes.

 To end this section, we mention the Dirichlet form concerning birth and death processes.
 Let
 $$
D(f)=\sum_{i=0}^{\ift}\mu_ib_i(f_i-f_{i+1})^2,
$$
and $\D^{\max}(D)=\set{f\in L^2(\mu): D(f)<\ift}$. Then it is proven in \ct[Proposition 1.3]{cbk3} that $(D,\D^{\max}(D))$ is regular if and only if
\be\lb{u1}
\sum_{i=0}^\ift \[\fr1{\mu_ib_i}+\mu_i\]=\ift.
\de
In other words, the Dirichlet form corresponding to $Q$ is unique iff \rf{u1} holds.

We remark that when $\ift$ is the exit or entrance boundary, the Dirichlet form is unique. Indeed, from \ct[Section 8.1]{and}, we know
that the equivalence condition for the exit boundary is $R<\ift, \mu=\ift, \sum_{i=0}^\ift 1/\mu_ib_i<\ift$;
 the equivalence condition for the entrance boundary is $S<\ift, \mu<\ift, \sum_{i=0}^\ift 1/\mu_ib_i=\ift$. Thus in any case, \rf{u1} holds.
For regular boundary ($R<\ift,S<\ift$), the problem is that the
Dirichlet form is not unique. We need to develop new technique other
than that used in this paper. For the natural boundary
($R<\ift,S=\ift$), the situation is different. Although the
Dirichlet form is unique, we will face the difficulty of the
essential spectrum problem. So the formula must be totally different
from that in this paper.

The rest of the paper is organized as follows. In Section 2, we derive the distributions of the hitting times
$T_{i,n}$ for finite birth and death process from Theorem \ref{t01}. In Section 3, we give the distribution of the
life time for the birth and death processes with $\ift$ the exit boundary, starting from any $i\ge0$. In Section 4,
we give the distributions of $T_{i,n}(i\ge n)$ for the birth and death processes with $\ift$ the entrance boundary.
And finally in Section 5, the distribution of the (fastest) strong stationary time is derived and we also study the
convergence in separation for the process.
\section{The finite state space}

Let's first solve Case I from Theorem \ref{t01}.

\bg{cor}\lb{c01}
For $0\le i<n<\ift$,
\be\lb{x1}
\E e^{-sT_{i,n}}=\fr{\prodd_{\nu=1}^n\fr{\lmd_\nu^{(n)}}{s+\lmd_\nu^{(n)}}}{\prodd_{\nu=1}^i\fr{\lmd_\nu^{(i)}}{s+\lmd_\nu^{(i)}}}, s\ge0.
\de
In particular,
\be\lb{r3}
\E T_{0,n}=\sum_{1\le\nu< n}\fr1{\lmd_\nu^{(n)}}, \q \E T_{i,n}=\sum_{1\le\nu< n}\fr1{\lmd_\nu^{(n)}}-\sum_{1\le\nu< i}\fr1{\lmd_\nu^{(i)}}.
\de
\end{cor}

\bg{proof}
Since $T_{0,n}=T_{0,i}+T_{i,n}$ by the property of the birth and death process and  $T_{0,i}, T_{i,n}$ are
independent by the strong Markov property, the corollary follows immediately from Theorem \ref{t01}.

\rf{r3} follows from \rf{x1} by a standard method to derive the moments from the Laplace transform.
\end{proof}

We remark that \rf{r3} can be called the eigentime identity for absorbed birth and death processes. Cf. \ct{m06}.
It is different from that in \ct[Chapter 3]{af}, where the eigentime identity for the ergodic finite Markov chain is involving in
the average hitting time. In \ct{m04}, this kind of eigentime identity for the continuous-time Markov chain on countable state space was studied.
See Section 5 below.

For $0\le n<N<\ift$, let $\hat\lmd_{n,1}^{(N)}<\hat\lmd_{n,2}^{(N)}< \cdots<\hat\lmd_{n,N-n}^{(N)}$ be the positive eigenvalues of $-\hat Q_n^{(N)}$, 
where
\be\lb{qnn}
\hat Q_n^{(N)}:=\left(
  \begin{array}{ccccccc}
    -({ a_{n+1}}+b_{n+1}) & b_{n+1} & 0& 0 & \cdots & 0& 0 \\
     a_{n+2} & -(a_{n+2}+b_{n+2}) & b_{n+2} & 0& \cdots & 0& 0 \\
    \vdots & \ddots & \ddots & \ddots & \ddots & \ddots&  \vdots\\
    0 & 0 & 0& 0 & \cdots &  a_{N} & -a_{N} \\
  \end{array}
\right).
\de
This is the generator of the birth and death process on $\set{n,\cdots, N}$ with absorbing state $n$ and
 reflecting state $N$. 
 
 As we mentioned in the last section, by taking the mapping on the state
space: $j\rar j': j'=N-j$, we can easily obtain the following results from Theorem \ref{t01} and Corollary \ref{c01}.

\bg{cor}\lb{c02}
For $0\le n<i\le N<\ift$,
$$
\E e^{-sT_{i,n}}=
\fr{\prodd_{\nu=1}^{N-n}\fr{\hat\lmd_{n,\nu}^{(N)}}{s+\hat\lmd_{n,\nu}^{(N)}}}{\prodd_{\nu=1}^{N-i}\fr{\hat\lmd_{i,\nu}^{(N)}}{s+\hat\lmd_{i,\nu}^{(N)}}}, s\ge0.
$$
In particular
$$
\E e^{-sT_{N,n}}=\prodd_{\nu=1}^{N-n}\fr{\hat\lmd_{n,\nu}^{(N)}}{s+\hat\lmd_{n,\nu}^{(N)}}, s\ge0
$$
and
\be\lb{z1}
\sum_{\nu=1}^{N-n}\fr1{\hat\lmd_{n,\nu}^{(N)}}=\E T_{N,n}=\sum_{j=n}^N\fr{1}{\pi_jb_j}\sum_{i=j+1}^N\pi_i.
\de
\end{cor}
\bg{proof}The second equality in \rf{z1} can be found on \ct[Page 264]{and}.
\end{proof}

\section{$\ift$ is the exit boundary}

In this section, we will derive the life time distribution for the minimal birth and death process when $Q$-processes are not unique. Let's
recall the facts about the uniqueness of birth and death process, see for example \ct{and,cbk1}.
And then we will study
the spectral theory for the minimal birth and death processes. The spectral theory helps
us pass from finite states to infinite states, especially we will establish what are the limits of the eigenvalues for finite birth and death processes
when the states go up to the infinity.


When $R<\ift$, the corresponding $Q$-processes
are not unique, for details see \ct[Chapter 8]{and} or \ct[Chapter 4]{cbk1}.
 Let $(X_t,t\ge0)$ be the corresponding continuous-time Markov chain
 with the minimal $Q$-function $P(t)=(p_{ij}(t):i,j\in E)$,
 that is,
 $$
 p_{ij}(t)=\P_i[X_t=j,t<\zeta]
 $$
 with $\zeta=\lim_{n\rar\ift}\xi_n$ the life time, where $\xi_n$ be the successive jumps:
$$
\xi_0=0,\q \xi_n=\inf\set{t: t>\xi_{n-1}, X_t\not=X_{\xi^{}_{n-1}}},\q n\ge1.
$$

\bg{prop}\lb{pr1}
\begin{description}
\item[(i)~] When staring from the state $0$, the life time $\zeta=\lim_{n\rar\ift}T_{0,n}$ a.s.


 \item[(ii)~] $\P_i[\zeta=\ift]=1$ or $0$ for any $i\in E$ and  $\E_0 \zeta=R$.
\end{description}
\end{prop}
\bg{proof} (i)~Let $X_0=0$. Note that $\zeta$ is the (first) time that the process jumps infinite times,
for any $n$, before $T_{0,n}$ the process jumps only finite times, so that $\zeta\ge T_{0,n}$ for any $n$. Thus
$\zeta\ge\lim_{n\rar\ift}T_{0,n}$ a.s. Conversely, since the birth and death process jumps once to two nearest neighbors, then $\xi_n\le T_{0,n}$ for any $n$. Thus
$\zeta=\lim_{n\rar\ift}\xi_n\le\lim_{n\rar\ift}T_{0,n}$ a.s.

(ii)~ See \ct[Chapter 8]{and} or \ct[Chapter 4]{cbk1}.
\end{proof}

Denote by $L^2(\mu)$ the usual (real) Hilbert space on $E=\set{0,1,2,\cdots}$.
Then it is well known that $Q^{(n)}, Q, P(t)$ are self-adjoint operators on $L^2(\mu)$.
For a self-adjoint operator $A$ on $L^2(\mu)$, denote $\sgm(A), \sgm_\ess(A)$ respectively the
spectrum and the essential spectrum of $A$. Here, the essential spectrum consists of continuous spectrum and
eigenvalues with infinite multiplicity.
When $\sgm_\ess(Q)=\emptyset$, denote by $\lmd_1<\lmd_2<\cdots$
all the eigenvalues of $-Q$. Actually, all the eigenvalues under consideration in this paper are of one multiplicity, 
for this see Theorem \ref{t24} below.

 The following result in \ct{m06} is our start point of the spectral theory for the minimal birth and death process.

 \bg{thm}\lb{m061} If $R<\ift, S=\ift$ (the exit boundary),
 then $P(t)$ is  a Hilbert-Schmidt operator for any $t>0$. So that $\sgm_\ess(Q)=\emptyset$ and
$$
\sum_{n\ge1}\lmd_n^{-1}=R.
$$
\end{thm}

To get the distribution of the life time, we need the following minimax principle for eigenvalues,
which is a variant of classical Courant-Fischer theorem for symmetric matrices. See  for example \ct[p.149]{hj}.

\bg{prop}\lb{pr2} Assume that $\sgm_\ess(Q)=\emptyset$.
Let $\lmd_\nu(\nu\ge1), \lmd_\nu^{(n)} (1\le\nu< n)$ be eigenvalues for $-Q^{(n)}$ in \rf{eq0} and $-Q$ in \rf{eq1} respectively.
Then for $ \nu\ge1$
\be
\lmd_\nu =\max_{f_1,\cdots,f_{\nu-1}\in L^2(\mu)}\min\set{D(f): \mu(f^2)=1, \mu(ff_i)=0,1\le i<\nu}.
\de
and
for $1\le \nu< n$
\be
\lmd_\nu^{(n)}=\max_{f_1,\cdots,f_{\nu-1}\in L^2(\mu)}\min\set{D(f):f|_{[n,\ift)}=0, \mu(f^2)=1, \mu(ff_i)=0,1\le i<\nu},
\de
\end{prop}

\bg{proof}
a)~ We prove the assertion for $\lmd_\nu$ first. Let $e_\nu$ be the corresponding eigenfunction for $\lmd_\nu$, then
$$
\lmd_\nu=\inf\set{D(f):\mu(f^2)=1,\mu(fe_j)=1,1\le j<\nu}.
$$
Since $f=\sum_{j=1}^\ift\mu(fe_j)e_j$ and $-Qe_j=\lmd_j e_j$, then
$$
D(f)=\mu((-Qf)f)=\sum_{j=1}^\ift\lmd_j\mu(fe_j)^2,
$$
so that for $f_i\in L^2(\mu)(1\le i<\nu)$,
$$\aligned
&\inf\set{D(f):\mu(f^2)=1,\mu(ff_i)=0,1\le i<\nu}\\
&=\inf\set{\sum_{j=1}^\ift\lmd_j\mu(fe_j)^2:\mu(f^2)=1,\mu(ff_i)=0,1\le i<\nu}\\
&\le\inf\set{\sum_{j=1}^\nu\lmd_j\lan f,e_j\ran^2:\mu(f^2)=1,\mu(ff_i)=0 (1\le i<\nu), \mu(fe_j)=0 (j>\nu)}\\
&\le\sup\set{\sum_{j=1}^\nu\lmd_j\lan f,e_j\ran^2:\mu(f^2)=1, \mu(fe_j)=0,j>\nu}\\
&=\lmd_\nu.
\endaligned
$$
Therefore,
$$
\lmd_\nu\ge\max_{f_1,\cdots,f_{\nu-1}\in L^2(\mu)}\min\set{D(f):\mu(f^2)=1, \mu(ff_i)=0,1\le i<\nu}.
$$
If we choose $f_i=e_i,1\le i<\nu$, the above equality holds.

b)~For $f|_{[n,\ift)}=0$,
\be
D(f)=\sum_{i=0}^{n-2}\mu_ib_i(f_i-f_{i+1})^2+\mu_{n-1}b_{n-1}f_{n-1}^2=:D^{(n)}(f),
\de
then it's easy to check that
 $D^{(n)}(f)$ is Dirichlet form for $Q^{(n)}$ in \rf{eq1}. The rest of the proof is the same as above.
\end{proof}

 Let
\be\lb{kei}
\scr{K}=\set{f:f \text{~has finite support}}.
\de
Define
$$
D(f)=\sum_{i=0}^{\ift}\mu_ib_i(f_i-f_{i+1})^2
$$
with the minimal domain $\scr{D}(D)$ consisting of the functions in the closure of $\K$ with respect to the norm
$||\cdot||_D: ||f||_D=\mu(f^2)+D(f)$. In this paper, we deal with the minimal Dirichlet form or the minimal processes, cf. \ct[Proposition 6.59]{cbk1}.

 \bg{thm}\lb{t24} Assume that $R<\ift, S=\ift$, then
$$\fa \nu\ge1, \lmd_\nu^{(n)}\downarrow\lmd_\nu.$$
Moreover, all eigenvalues $\lmd_\nu$ are distinct (each of one multiplicity).
\end{thm}

\bg{proof} a)~ We use Proposition \ref{pr2} to prove the monotonicity for $\lmd_\nu^{(n)}$.
For any fixed $1\le\nu<n$ and $f_1,\cdots,f_{\nu-1}\in L^2(\mu)$, if $f$ is such that $f|_{[n,\ift)}=0,\mu(f^2)=1, \mu(ff_i)=0,1\le i<\nu$, then
 $f|_{[n+1,\ift)}=0,\mu(f^2)=1, \mu(ff_i)=0,1\le i<\nu$, so that
 $$\aligned
 &\min\set{D(f):f|_{[n,\ift)}=0, \mu(f^2)=1, \mu(ff_i)=0,1\le i<\nu}\\
 &\ge \min\set{D(f):f|_{[n+1,\ift)}=0, \mu(f^2)=1, \mu(ff_i)=0,1\le i<\nu}
 \endaligned
 $$
 and $\lmd_\nu^{(n)}\ge\lmd_\nu^{(n+1)}$ for  $1\le\nu<n$. This proves the monotonicity. Thus the limit $\lim_{n\rar\ift}\lmd_\nu^{(n)}=:\hat\lmd_\nu$
 exists for any $\nu\ge1$.

b)~As pointed in Section 1, when $R<\ift, \mu=\ift$, $(D,\D(D))$ is a regular Dirichlet form. That is, let
$$
\D_{\max} (D)=\set{f\in L^2(\mu): D(f)<\ift},
$$
and $(D,\D_{\max}(D))$ be the maximum Dirichlet form, then $\D(D)=\D_{\max}(D)$. Therefore, when $R<\ift, \mu=\ift$ it
follows from Proposition \ref{pr2} that
\be\lb{l01}
\lmd_\nu =\max_{f_1,\cdots,f_{\nu-1}\in L^2(\mu)}\min\set{D(f): f\in \K, \mu(f^2)=1, \mu(ff_i)=0,1\le i<\nu}.
\de

c)~On one hand, we have $\lmd_\nu\le\hat\lmd_\nu$. Indeed, it follows from  \rf{l01} and monotonicity that for $\nu\ge1$,
 $$\aligned
 \lmd_\nu&=\max_{f_1,\cdots,f_{\nu-1}\in L^2(\mu)}\min\set{D(f):f\in \K, \mu(f^2)=1, \mu(ff_i)=0,1\le i<\nu}\\
 &=\max_{f_1,\cdots,f_{\nu-1}\in L^2(\mu)}\min\set{D(f):\ex n \ge0, f|_{[n,\ift)}=0; \mu(f^2)=1, \mu(ff_i)=0,1\le i<\nu}\\
 &=\max_{f_1,\cdots,f_{\nu-1}\in L^2(\mu)}\inf_{n\ge1}\min\set{D(f):f|_{[n,\ift)}=0, \mu(f^2)=1, \mu(ff_i)=0,1\le i<\nu}\\
 &\le\inf_{n\ge1}\max_{f_1,\cdots,f_{\nu-1}\in L^2(\mu)}\min\set{D(f):f|_{[n,\ift)}=0, \mu(f^2)=1, \mu(ff_i)=0,1\le i<\nu}\\
 &=\lim_{n\rar\ift}\max_{f_1,\cdots,f_{\nu-1}\in L^2(\mu)}\min\set{D(f):f|_{[n,\ift)}=0, \mu(f^2)=1, \mu(ff_i)=0,1\le i<\nu}\\
 &=\lim_{n\rar\ift}\lmd_\nu^{(n)}=\hat\lmd_\nu.
 \endaligned
 $$
 On the other hand, from Corollary \ref{r3} we have
 $$
 \E T_{0,n}=\sum_{1\le\nu< n}\fr1{\lmd_\nu^{(n)}}=\sum_{1\le\nu<\ift}\fr1{\lmd_\nu^{(n)}}I_{[\nu<n]},
 $$
 then it follows from monotone convergence theorem that
 \be\lb{r1}
 R=\E_0 \zeta=\sum_{1\le\nu<\ift}\fr1{\hat\lmd_\nu}.
 \de
 But  we already know from Theorem \ref{m061} that
 \be\lb{r2}
 R=\sum_{1\le\nu<\ift}\fr1{\lmd_\nu}<\ift.
 \de
 Since $\lmd_\nu\le \hat\lmd_\nu (\nu\ge1)$, it must hold that $\lmd_\nu= \hat\lmd_\nu$ for any $\nu\ge1$.

d)~Next we will prove that all eigenvalues $\set{\lmd_\nu,\nu\ge1}$ are distinct. For this we only need to prove that the eigenspace for any $\lmd_\nu$ is
of one dimension.
 Indeed, let $-\lmd$ be an eigenvalue and $g$ the corresponding eigenfunction.
From $Qg(i)=-\lmd g_i, i\ge0$, we have
\be\lb{eq11}
b_0(g_1-g_0)=-\lmd g_0, a_{i}(g_{i-1}-g_i)+b_i(g_{i+1}-g_i)=-\lmd g_i, i\ge1.
\de
Since $\mu_ib_i=\mu_{i+1}a_{i+1} (i\ge0)$, it follows from \rf{eq11} that
$$
g_{k+1}=-\fr{\lmd}{\mu_kb_k}\sum_{i=0}^k\mu_ig_i+g_k, k\ge0.
$$
(cf.\ct{cbk3}) This means that eigenfunction $g$ is determined uniquely once $g_0$ is given.

\end{proof}

\bg{rem}
For the first eigenvalue ($\nu=1$), it was proved by Chen M.-F.(2009) without the assumption $\sgm_\ess(Q)=\emptyset$
in case that $\lmd_1$ is defined by the classical Poincar\'e variational formula:
$$
\lmd_1=\inf\set{ \sum_{i=0}^\ift\mu_ib_i(f_i-f_{i+1})^2: f\in\K, \sum_{i=0}^\ift\mu_if_i^2=1}.
$$
\end{rem}

\bg{thm}\lb{t36} Assume $R<\ift,S=\ift$. Let $\zeta$ be the life time for the minimal process, then
$$
\E_0 e^{-s\zeta}=\prod_{\nu=1}^\ift\fr{\lmd_\nu }{s+\lmd_\nu }, s\ge0.
$$
And for any $i\ge0$, let $T_{i,\ift}=\lim_{n\rar\ift}T_{i,n}$, then
$$
\E_i e^{-s\zeta}=\E e^{-sT_{i,\ift}}=\fr{\prodd_{\nu=1}^\ift\fr{\lmd_\nu }{s+\lmd_\nu }}{\prodd_{\nu=1}^i\fr{\lmd_\nu^{(i)}}{s+\lmd_\nu^{(i)}}}, s\ge0.
$$
\end{thm}
\bg{proof}
The assertions follow from the monotone convergence theorem and Proposition \ref{pr1}, Corollary \ref{c01}, Theorem \ref{t24}.
\end{proof}

For the exit boundary, the distributions of $T_{i,n}$ for $0\le i,n\le\ift$ are all known. When $0\le i< n<\ift$, the
distribution is given by Corollary \ref{c01}, while $0\le i< n=\ift$, the
distribution is given by Theorem \ref{t36}. 

\section{$\ift$ is the entrance boundary}

In this section we will deal with Case IV for the birth and death process with $\ift$ the entrance boundary, i.e. $R=\ift, S<\ift$, and the
corresponding  $Q$-process
is unique.

Since $S<\ift$, $\mu<\ift$.
Let $\pi_i=\mu_i/\mu$, then $\pi=(\pi_i,i\ge0)$ is a probability measure on $E$, so that the process is reversible with
respect to $\pi$. Now we will consider the spectral theory for
 operators on Hilbert space $L^2(\pi)$.


 For $n\ge0$, let
\be\lb{qn}
\hat Q_n=\left(
  \begin{array}{cccccc}
     -({a_{n+1}}+b_{n+1}) & b_{n+1} & 0 & \cdots & \cdots& \cdots \\
     a_{n+2} & -(a_{n+2}+b_{n+2}) & b_{n+2} & \cdots & \cdots& \cdots \\
     \ddots & \ddots & \ddots & \ddots & \ddots&  \ddots\\
  \end{array}
\right)
\de
be the generator of the birth and death process absorbed at state $n$.

Let $\hat\pi^{(n)}=(\pi_i:i>n)$ and $\hat E_n=\set{n+1,n+2,\cdots}$.
It is easy to check that $\hat Q_n$ is symmetric with respect to $\hat \pi{(n)}$ and then $\hat Q_n$ is a self-adjoint operator
in $L^2(\hat E_n,\hat \pi^{(n)})$.
When $\sgm_\ess(\hat Q_n)=\emptyset$, denote by $\hat\lmd_{n,1}<\hat\lmd_{n,2}<\cdots$
all the positive eigenvalues of $-\hat Q_n$, as we know from Theorem \ref{t24} that each eigenvalue is of one multiplicity.
When $n=0$, the subscript $0$ is dropped.

\bg{thm}\lb{t31} For the birth and death process with $\ift$ the entrance boundary, i.e. $R=\ift, S<\ift$, then
   $\sgm_\ess(\hat Q_n)=\emptyset$, and  for any $n\ge0$
\be\lb{x2}
  S_n:=\sum_{j=n}^\ift\fr{1}{\pi_jb_j}\sum_{i=j+1}^\ift\pi_i=\sum_{\nu\ge1}\hat\lmd_{n,\nu}^{-1}<\ift.
\de
\end{thm}
\bg{proof} It follows from \ct[Theorem 1.4]{m04} that
$
\sgm_\ess(Q)=\emptyset\q \text{and}\q \sum_{\nu\ge1}\hat\lmd_{\nu}^{-1}<\ift,
$
where $Q$ is defined by \rf{eq0} and $\set{\lmd_\nu:\nu\ge0}$ is spectrum of $-Q$ in $L^2(\pi)$ with $\lmd_0=0$. But
since $Q$ and $\hat Q_n$ differ only from a finite states, their essential spectrum is same (se for example \ct[Theorem 5.35 on page 244]{kato}).

Now we prove the identity in \rf{x2}. Let for $i,j>n$
$$
p^{(n)}_{ij}(t)=\P_i[X_t=j,t<T_{i,n}], g^{(n)}_{ij}=\int_0^\ift p^{(n)}_{ij}(t)dt.
$$
By a similar method as in the proof of \ct[Theorem 1.4]{m06}, we can get that
$$
\sum_{\nu\ge1}\hat\lmd_{n,\nu}^{-1}=\sum_{i>n}g^{(n)}_{ii}=\sum_{i>n}\fr{1}{(a_i+b_i)\P_i[\tau_i^+=\ift]}.
$$
Here  $\tau_j^+=\inf\set{t\ge\text{the first jump time}:X_t=j}$ is the return time. 
Note that once $ [\tau_i^+=\ift|X_0=i]$ happens, it must first jump to state $i-1$,
otherwise it can be back to state $i$ in finite time almost surely since the original $Q$-process is ergodic. 
Next when it comes to state $i-1$, it must arrive to state $n$ before it arrives state $i$. Thus
we have
$$
\P_i[\tau_i^+=\ift]=\fr{a_i}{a_i+b_i}\P_{i-1}[T_{i-1,n}<\tau_i^+],
$$
and by a standard martingale method and letting $s_i=\sum_{j<i}(\pi_jb_j)^{-1}$ be the scale function (see for example \ct[Theorem 4 in \S 7.11]{wy}), 
we obtain
$$
\P_{i-1}[T_{i-1,n}<\tau_i^+]=\fr{s_i-s_{i-1}}{s_i-s_n}=\fr{(\mu_{i-1}b_{i-1})^{-1}}{\sum_{n\le j<i}(\pi_jb_j)^{-1}}
=\fr{(\mu_{i}a_{i})^{-1}}{\sum_{n\le j<i}(\pi_jb_j)^{-1}}.
$$
Therefore
$$
\sum_{\nu\ge1}\hat\lmd_{n,\nu}^{-1}=\sum_{i>n}\pi_i\sum_{n\le j<i}\fr{1}{\pi_jb_j}=\sum_{j=n}^\ift\fr{1}{\pi_jb_j}\sum_{i=j+1}^\ift\pi_i.
$$
\end{proof}

As in the last section, we also need the minimax principle for eigenvalues of $\hat Q_n$.

\bg{prop}\lb{pr3}
For $0\le n<N<\ift$, let
$$
D_n^{(N)}(f)=\mu(-f\hat Q_n^{(N)}f)=\sum_{i=n+1}^{N-1}\mu_ib_i(f_i-f_{i+1})^2+\mu_nb_nf_{n+1}^2.
$$
Then for $1\le \nu< N-n$
\be
\hat\lmd_{n,\nu}^{(N)}=\max_{f_1,\cdots,f_{\nu-1}\in L^2(\mu)}\min\set{D_n^{(N)}(f): \mu(f^2)=1, \mu(ff_i)=0,1\le i<\nu}.
\de
\end{prop}
\bg{proof}
The proof is direct and is omitted.
\end{proof}

\bg{prop}\lb{pr4} Assume that $\sgm_\ess(Q)=\emptyset$.
Let $\hat\lmd_{n,\nu}$ be eigenvalues for $-\hat Q^{(n)}$ in \rf{qn}.
Then for $ \nu\ge1$
\be\lb{44}
\hat\lmd_{n,\nu}=\max_{f_1,\cdots,f_{\nu-1}\in L^2(\mu)}\min\set{D(f): f_{[0,n]}=0,\mu(f^2)=1, \mu(ff_i)=0,1\le i<\nu}.
\de
\end{prop}
\bg{proof}
Since $D(f)=\sum_{i=n+1}^{\ift}\mu_ib_i(f_i-f_{i+1})^2+\mu_nb_nf_{n+1}^2=\mu(-f\hat Q_nf)$. The rest of proof is similar
to that of Proposition \ref{pr2}
\end{proof}

We need the approximation procedure when $N\rar\ift$. For this purpose, we need to do more.

Fix $n\ge0$ and define
\be\lb{kap}
\hat{ \scr{K}}=\set{f\in L^\ift: \set{f\not=0}\subset \set{n+1,\cdots,N}\text{~ for some~$N$}}.
\de
and $\hat {\scr {K}}_{L}=\set{g:=cf+d: f\in\hat{\scr{K}},c,d\in \R}$.
Define
$$
D(f)=\sum_{i=0}^{\ift}\mu_ib_i(f_i-f_{i+1})^2
$$
with the domain $\scr{D}(D)$ consisting of the functions in the closure of $\hat\K_L$ with respect to the norm
$||\cdot||_D: ||f||_D=\mu(f^2)+D(f)$.

Since $\ift$ is the entrance boundary, the Dirichlet form is unique as explained in Section 1, thus
$\scr{D}(D)=\D^{\max}(D)=\set{f\in L^2(\mu): D(f)<\ift}$ and we can rewrite \rf{44} as
$$
\hat\lmd_{n,\nu}=\max_{f_1,\cdots,f_{\nu-1}\in L^2(\mu)}\min\set{D(f): f\in\hat\K_L,\mu(f^2)=1, \mu(ff_i)=0,1\le i<\nu}.
$$
This leads to

 \bg{thm}\lb{zz}Assume that $R=\ift, S<\ift$, then
$$\fa \nu\ge1, \hat\lmd_{n,\nu}^{(N)}\downarrow\hat\lmd_{n,\nu} \q \text{as} \q N\rar\ift.
$$
\end{thm}
\bg{proof}(a)~Let $f\in \hat\K_L$, assume that $f_{[N,\ift)}=c$, then
$$\aligned
D(f)&=\sum_{i=0}^{\ift}\mu_ib_i(f_i-f_{i+1})^2=\sum_{i=n+1}^{\ift}\mu_ib_i(f_i-f_{i+1})^2+\mu_nb_nf_{n+1}^2\\
   &=\sum_{i=n+1}^{N-1}\mu_ib_i(f_i-f_{i+1})^2+\mu_nb_nf_{n+1}^2=D_n^{(N)}(f).
   \endaligned
$$
Here in $D_n^{(N)}(f)$, $f$ is viewed as a function on $\set{n+1,\cdots,N}$. From this, we can easily deduce the monotonicity of
$\hat\lmd_{n,\nu}^{(N)}$ in $N$.

(b)~By Theorem 9.11 in \ct{cbk1} and Proposition \ref{pr4} above, we have that for $\nu\ge1$,
 $$\aligned
 \hat\lmd_{n,\nu}&=\max_{f_1,\cdots,f_{\nu-1}\in L^2(\mu)}\min\set{D(f):f\in \hat\K_L, \mu(f^2)=1, \mu(ff_i)=0,1\le i<\nu}\\
 &=\max_{f_1,\cdots,f_{\nu-1}\in L^2(\mu)}\min\left\{D_n^{(N)}(f):\ex N>n, f|_{[N,\ift)}=\text{constant}; \mu(f^2)=1,\mu(ff_i)=0,1\le i<\nu\right\}\\
 &\le\lim_{N\rar\ift}\max_{f_1,\cdots,f_{\nu-1}\in L^2(\mu)}\min\set{D_n^{(N)}(f):f|_{[N,\ift)}=\text{constant}, \mu(f^2)=1, \mu(ff_i)=0,1\le i<\nu}\\
 &=\lim_{N\rar\ift}\hat\lmd_{n,\nu}^{(N)}:=\alp_{n,\nu}.
 \endaligned
 $$

 (c)~It follows from \rf{z1} in Corollary \ref{c02} that
 $$
 \sum_{j=n}^N\fr{1}{\pi_jb_j}\sum_{i=j+1}^N\pi_i=\sum_{\nu\ge1}\fr{1}{\hat\lmd_{n,\nu}^{(N)}},
 $$
from which by letting $N\rar\ift$ we have
$$
\sum_{j=n}^\ift\fr{1}{\pi_jb_j}\sum_{i=j+1}^\ift\pi_i=\sum_{\nu\ge1}\fr{1}{\alp_{n,\nu}}.
$$
But for any $\nu\ge1$, $\hat\lmd_{n,\nu}\le\alp_{n,\nu}$ and it follows from Theorem \ref{t31} that
$$
\sum_{j=n}^\ift\fr{1}{\pi_jb_j}\sum_{i=j+1}^\ift\pi_i=\sum_{\nu\ge1}\fr{1}{\hat\lmd_{n,\nu}}.
$$
Thus it must hold that $\hat\lmd_{n,\nu}=\alp_{n,\nu}=\lim_{N\rar\ift}\hat\lmd_{n,\nu}^{(N)}$ for any $\nu\ge1$.

\end{proof}

\bg{thm}\lb{t45} Assume the birth and death process is such that $\ift$ the entrance boundary. For $n\ge1$, we have
\be\lb{my}
\E e^{-sT_{n,0}}=
\fr{\prodd_{\nu=1}^\ift \fr{\hat\lmd_\nu}{s+\hat\lmd_{\nu}}}{\prodd_{\nu=1}^\ift \fr{\hat\lmd_{n,\nu}}{s+\hat\lmd_{n,\nu}}},s\ge0.
\de
Let $T_{\ift,0}=\lim_{n\rar\ift}T_{n,0}$, then
\be\lb{mz}
\E e^{-sT_{\ift,0}}=
 {\prodd_{\nu=1}^\ift \fr{\hat\lmd_\nu}{s+\hat\lmd_{\nu}}},s\ge0.
\de
\end{thm}

\bg{proof}
a)~By using the monotone convergence theorem, we can get \rf{my} from Theorem \ref{zz} and Corollary \ref{c02} immediately.

b)~ To pass from \rf{my} to \rf{mz}, we only need to show that $\lim_{n\rar\ift}\hat\lmd_{n,1}=\ift$, so that $\hat\lmd_{n,\nu}$ tend to infinity
uniformly in $\nu\ge1$ as $n\rar\ift$. This is a consequence of \ct[Corollary 1.1]{m04} and \ct[Theorem 3.4]{m061}.
\end{proof}

For the entrance boundary, the distributions of $T_{i,n}$ for $0\le i,n\le\ift$ are all known. When $0\le i< n<\ift$, the
distribution of $T_{i,n}$ is given by Corollary \ref{c01}, while $0\le i< n=\ift$, $T_{i,\ift}=\ift$ a.s. since the life time
$\zeta=\ift$ a.s.
When $0\le n< i\le\ift$, the distribution of $T_{i,n}$ is given by Theorem \ref{t45}.

\section{Application to the fastest strong stationary time}

In this section, we apply the theorems to the strong stationary time. We give the distribution to the fastest
strong stationary times, and then study the exponential convergence in separation for the birth and death process.

First of all, we would like to recall some facts about the strong ergodicity for the birth and death process.

 Let $T$ be the average hitting time:
$
T=\sum_{ij}\pi_i\pi_j \E T_{i,j}.
$
(cf.\ct[Chapter 3]{af}) From \ct[Chapter 8]{and}, we can eventually calculate out that (cf.\ct{m04})
$$
T=\sum_{k=0}^\ift\fr{1}{\pi_kb_k}\sum_{i=0}^k\pi_i\sum_{i=k+1}^\ift\pi_i.
$$

In the following theorem, we summary the facts of strong ergodicity.

\bg{thm}\lb{t51}Assume that the process is unique (i.e. $R=\ift$).
The following statements are equivalent.
\begin{description}
\item [(i)]$S<\ift$.
  \item [(ii)]The process is strongly ergodic.
  \item[(iii)] $T<\ift.$
  \item[(iv)] $\sgm_\ess(Q)=\emptyset$ and $\sum_{\nu\ge1}\lmd_{\nu}^{-1}<\ift.$
  \end{description}
  Furthermore, when $\sgm_\ess(Q)=\emptyset$, then
  $$
  \sum_{\nu\ge1}\lmd_\nu^{-1}=T.
  $$
\end{thm}
\bg{proof}
The equivalence of (i) and (ii) was proved in \ct{zlh} and the other assertions were prove in \ct{m04}.
\end{proof}

  Now we study the distribution of the strong stationary time for the strongly ergodic birth and death process.
A strong stationary time (SST) is a (minimal) randomized stopping
time $\tau$ for $X_t$ such that $X_\tau$ has the distribution $\pi$ and it independent of $\tau$.

With the aid of Theorem \ref{t36} and the duality established in \ct{fill1}, we can obtain the following result,
which extends  the result in \ct{fill2} to the denumerable case.

\bg{thm}\lb{sst} For a strongly ergodic birth and death process (i.e. $R=\ift,S<\ift$), the SST $\tau$ has distribution
$$
\E_0 e^{-s\tau}=\prod_{\nu=1}^\ift\fr{\lmd_\nu}{s+\lmd_\nu},
$$
where $\set{\lmd_\nu:\nu\ge1}$ are the positive eigenvalues of $-Q$ for the strongly ergodic birth and death process.
\end{thm}
\bg{proof}We follow the argument in \ct[Section 3.3]{fill1} with a minor modification. Let
$$
H_i=\sum_{j\le i}\pi_j,
$$
where $\pi_j=\mu_j/\mu$ is as before. Then $0<1/\mu=\pi_0\le H_i\le1$.

Define the dual $Q^*$-birth and death process with parameters $(a^*_i,b_i^*)$ given by
\be\lb{**}
a_i^*=\fr{H_{i-1}}{H_i}b_i,\q b_i^*=\fr{H_{i+1}}{H_i}a_{i+1}.
\de
Also
$$
\mu_i^*=\fr{b_0^*\cdots b_{i-1}^*}{a_1^*\cdots a_i^*}=\fr{b_0}{\mu_ib_i}\fr{H_i^2}{H_0^2}.
$$
Note that $H_i\ge H_0$, we have
$$
\mu^*=\sum_{i=0}^\ift\mu_i^*=\sum_{i=0}^\ift\fr{b_0}{\mu_ib_i}\fr{H_i^2}{H_0^2}\ge\sum_{i=0}^\ift\fr{b_0}{\mu_ib_i}\fr{H_i}{H_0}=
\sum_{i=0}^\ift\fr{b_0}{\mu_ib_i}\sum_{j\le i}\mu_j=b_0R=\ift,
$$
and as $1/\mu\le H_j\le1, \mu_ib_i=\mu_{i+1}a_{i+1}$,
$$\aligned
R^*&=\sum_{i=0}^\ift\fr{1}{\mu_i^*b^*_i}\sum_{j\le i}\mu^*_j=\sum_{i=0}^\ift\fr{\mu_{i+1}}{b_0H_iH_{i+1}}\sum_{j\le i}\fr{b_0}{\mu_jb_j}H_j^2\\
&=\sum_{j=0}^\ift\fr{H_j^2}{\mu_jb_j}\sum_{ i\ge j}\fr{\mu_{i+1}}{H_iH_{i+1}}\le\mu^{-2}S<\ift.
\endaligned
$$
Thus the $Q^*$-process (minimal process) is with $\ift$ the exit boundary and $\tau$ has the same distribution as the life time $\zeta^*$ (the time attaining $\ift$)
for $Q^*$-process. Applying Theorem \ref{t36}, we have
$$
\E_0 e^{-s\tau}=\E_0 e^{-s\zeta^*}=\prod_{\nu=1}^\ift\fr{\lmd_\nu^* }{s+\lmd_\nu^*},
$$
where $\set{\lmd_\nu^*,\nu\ge1}$ is spectrum of $-Q^*$. To complete the proof, it suffices to prove that $\lmd_\nu^*=\lmd_\nu$ for $\nu\ge1$.

Actually, let the link matrix $\Lmd=(\Lmd_{ij})$ as $\Lmd_{ij}=1_{[j\le i]}\pi_j/H_i$, then $\Lmd Q=Q^*\Lmd$. This means that $Q$ has the same spectrum
in $L^2(\mu)$ as $Q^*$ in $L^2(\mu^*)$. But as proved in Theorems \ref{m061} and \ref{t51},  the spectra are all eigenvalues, so that
 $\lmd_\nu^*=\lmd_\nu$ for $\nu\ge1$.
\end{proof}

Theorem \ref{sst} enables one to obtain the unform convergence in separation. Actually we will prove that the uniform convergence in separation
is equivalent to the strong ergodicity for the birth and death process.

Let $p_{ij}(t)$ be the transition function for the
birth and death process with stationary distribution $\pi$, define the separation:
$$
s_i(t)=\sup_{j}\(1-\fr{p_{ij}(t)}{\pi_j}\), \fa i\ge0,t\ge0.
$$
For the elementary properties of separation, see \ct{fill1}. For example, we have the following relation for total variance distance and separation:
\be\lb{var}
||p_{i\cdot}(t)-\pi||_\var:=\sum_j|p_{ij}(t)-\pi_j|\le s_i(t).
\de
The main theorem in \ct{fill1} (see also \ct[Proposition 1]{fill2}) says
\be\lb{sep}
s_i(t)\le\P_i[\tau>t], \fa 0\le t<\ift.
\de
This leads to study the distribution of $\tau$ starting from any state other than $0$.

\bg{prop}\lb{pr5}Let $\tau$ be the SST in Theorem \ref{sst} and $X_t$ be the $Q$-process and $X_t^*$ be the minimal $Q^*$-process
with $Q^*$ given by \rf{**}. Then for $t\ge0$
\be\lb{mm}
\P_0[\tau\le t]=\P_0[\zeta^*\le t],
\de
and for $i\ge1$
\be\lb{mm1}
\P_i[\tau\le t]=\fr{1}{\pi_i}\bigg\{H_i\P_i[\zeta^*\le t]-H_{i-1}\P_{i-1}[\zeta^*\le t]\bigg\}.
\de
Consequently, for $i\ge1, \E_i\tau\le\E_{i-1}\zeta^*\le\E_0\zeta^*$ and
\be\lb{mm0}
\sup_{i\ge0}\E_i\tau=\E_0\tau.
\de
\end{prop}

\bg{proof}
The proof of \rf{mm} can be found in \ct{fill1}. We will prove \rf{mm1}. Let $m$ and $m^*$ be the distribution of $X_0$ and
$X_0^*$ respectively, it follow from \ct[Proposition 4]{fill1} that if $m=m^*\Lmd$, then
\be\lb{mm2}
\P_m[\tau\le t]=\P_{m^*}[\zeta^*\le t].
\de
For any $i\ge1$, let $m^*=\dlt_i$ the Dirac measure, then $m_k=\fr{\pi_k}{H_i}1_{[0\le k\le i]}$ and equality \rf{mm2} implies that
$$
\sum_{0\le k\le i}\fr{\pi_k}{H_i}\P_k[\tau\le t]=\P_i[\zeta^*\le t]\q\text{or}\q \sum_{0\le k\le i}{\pi_k}{}\P_k[\tau\le t]=H_i\P_i[\zeta^*\le t].
$$
Thus $$
\P_i[\tau\le t]=\fr{1}{\pi_i}\bigg\{H_i\P_i[\zeta^*\le t]-H_{i-1}\P_{i-1}[\zeta^*\le t]\bigg\}.
$$
Since $\P_i[\zeta^*\le t]$ increases in $i\ge0$, we have
$$
\P_i[\tau\le t]\ge\fr{1}{\pi_i}\bigg\{H_i\P_{i-1}[\zeta^*\le t]-H_{i-1}\P_{i-1}[\zeta^*\le t]\bigg\}=\P_{i-1}[\zeta^*\le t]\ge\P_{0}[\zeta^*\le t].
$$
Thus
$$
\P_i[\tau> t]\le\P_{i-1}[\zeta^*> t]\le \P_0[\zeta^*> t],
$$
which implies
$$
\E_i\tau=\int_0^\ift\P_i[\tau> t]dt\le\int_0^\ift\P_{i-1}[\zeta^*> t]dt=\E_{i-1}\zeta^*\le \E_0\zeta^*.
$$
\end{proof}

\bg{cor}\lb{k1}
We have for any $\ell\ge0$
\be\lb{k2}
\E_0\tau^\ell\le (\E_0\tau)^\ell/\ell!
\de
and for any $\lmd<1/\E_0\tau$
\be\lb{k3}
\E_0e^{\lmd\tau}\le(1-\lmd\E_0\tau)^{-1}.
\de
\end{cor}

\bg{proof}Since from \rf{mm1}, $\tau$ and $\zeta^*$ have the same distribution under $\E_0$, we need only prove \rf{k2} with $\tau$ replaced by $\zeta^*$.
It follows from \ct[Lemma 4.1]{m06} that  for any $\ell\ge0$
$$
\E_0[(\zeta^*)^\ell]\le (\E_0\zeta^*)^\ell/\ell!
$$
and for any $\lmd<1/\E_0\tau$
$$
\E_0e^{\lmd\tau}\le(1-\lmd\E_0\tau)^{-1}.
$$
\end{proof}

\bg{thm}\lb{sst1} Let
$$
S(t)=\sup_i s_i(t),
$$
and $\tau$ be the SST in Theorem \ref{sst}. Then the following statements are equivalent.
\begin{description}
  \item[(i)] The process is strongly ergodic.
  \item[(ii)] $\tau<\ift$ a.s.
  \item[(iii)] $\E_0 \tau<\ift$.
  \item[(iv)] $\lim_{t\rar\ift}S(t)=0$.
\end{description}
\end{thm}

\bg{proof}
The equivalence of (i)-(iii) follows from Theorems \ref{sst} and \ref{t31}, and Proposition \ref{pr5}.  We will prove the implication of
(iii)$\Rar$(iv) and (iv)$\Rar$(i).

(a)~Suppose that $\sup_i\E_i\tau<\ift$. By duality, we have $\E_0\zeta^*<\ift$, it follows from \rf{sep} and Markov inequality
that
$$
S(t) \le\sup_i\P_i[\tau>t]\le \fr{\sup_i\E_i\tau}{t}\le\fr{\E_0\tau}{t}\rar0, \q t\rar\ift.
$$

(b)~Suppose $\lim_{t\rar\ift}S(t)=0$. It follows from the inequality \rf{var} that
$$
\sup_i||p_{i\cdot}(t)-\pi||_\var\le S(t)\rar0, \q t\rar\ift,
$$
which is strong ergodicity.
\end{proof}

We remark that under any  assumption of $(1)-(4)$ above, letting
$$
\bt=\sup\set{\eps: \exists C<\ift, S(t)\le Ce^{-\eps t}, \fa t\ge0}
$$
be the optimal convergence rate in separation, then it follows easily from Corollary \ref{k1} and Theorem \ref{t51} that
$$
\bt\ge\fr {1}{\E_0\tau}=\(\sum_{\nu=1}^\ift\lmd_\nu^{-1}\)^{-1}=\(\sum_{k=0}^\ift\fr{1}{\pi_kb_k}\sum_{i=0}^k\pi_i\sum_{i=k+1}^\ift\pi_i\)^{-1}.
$$

{\bf Acknowledgement}\q The authors would thank Professors Mu-Fa Chen and Feng-Yu Wang for their valuable suggestions.


\end{document}